\documentclass[12pt,reqno]{amsart}
\usepackage{amsmath,amssymb,extarrows}
\usepackage{hyperref}
\usepackage{url}
\usepackage{tikz,enumerate}
\usepackage{diagbox}
\usepackage{appendix}
\usepackage{epic}
\usepackage{comment}
\usepackage{bm}
\usepackage{adjustbox}
\usepackage{booktabs}
\usepackage{array}

\usepackage{float} %support many floats

\usepackage{cite}
\usepackage{hyperref}
\usepackage{array}

\allowdisplaybreaks[4]

\newtheorem{theorem}{Theorem}

\newtheorem{lemma}[theorem]{Lemma}
\newtheorem{prop}[theorem]{Proposition}
\theoremstyle{definition}

\theoremstyle{remark}
\newtheorem{rem}{Remark}
\numberwithin{equation}{section}
\numberwithin{theorem}{section}
\numberwithin{defn}{section}
\newcommand{\e}{\mathrm{e}}

\begin{document}
\title[New lower bounds for even and odd values of the partition function]
 {New Lower Bounds for  Even and Odd Values of \\ the Partition Function}

\author{Liuquan Wang}
\address{School of Mathematics and Statistics, Wuhan University, Wuhan 430072, Hubei, People's Republic of China}
\email{wanglq@whu.edu.cn; mathlqwang@163.com}

\subjclass[2020]{05A17, 11P83, 11P82, 11F03, 11F30.}

\keywords{Partition function; parity; Dedekind eta function; Faber polynomials; partition congruences}

%\date{November 6, 2026}
%\dedicatory{}
%\thanks{}

\begin{abstract}
Let $N_{\mathrm{odd}}(X)$ and $N_{\mathrm{even}}(X)$ denote the numbers of integers $0\leq  n\leq X$ for which the ordinary partition function $p(n)$ is odd and even, respectively. The Parkin--Shanks conjecture predicts that $N_{\mathrm{odd}}(X)$ and $N_{\mathrm{even}}(X)$ are both asymptotic to $X/2$. Bella\"iche and Nicolas proved that $N_{\mathrm{even}}(X)\gg \sqrt{X} \log \log X$, while Zheng, building on work of Ono and Swaminathan, recently proved that $N_{\mathrm{odd}}(X)\gg \sqrt{X}$. We improve both bounds simultaneously by proving $\min\{N_{\mathrm{odd}}(X), N_{\mathrm{even}}(X)\}\gg X^{3/4}/(\log X \log \log X)$. These results have been verified formally in Lean by AxiomProver.
\end{abstract}

\maketitle

\section{Introduction}
For any positive integer $n$, a partition of $n$ is a nonincreasing sequence of positive integers whose sum is $n$.  We denote the number of partitions of $n$ by $p(n)$ with the convention that $p(0)=1$. For example, we have $p(4)=5$ as there are five partitions of 4:
$$4, \quad 3+1, \quad 2+2, \quad 2+1+1, \quad 1+1+1+1.$$
Euler's generating function identity states that
\begin{align}
 \sum_{n\geq 0}p(n)q^n=\prod_{n \geq 1}\frac{1}{1-q^n}, \quad |q|<1.
\end{align}

In 1919, Ramanujan \cite{Ramanujan} discovered the celebrated congruences: for all $n\geq 0$,
\begin{align}\label{pn-cong}
    p(5n+4)\equiv 0 \!\!\!\pmod{5}, ~ p(7n+5)\equiv 0 \!\!\! \pmod{7}, ~p(11n+6)\equiv 0  \!\!\!\pmod{11}.
\end{align}
He \cite{Ramanujan1920} further speculated that no other congruences of this simple form exist.  Ahlgren and Boylan  \cite{Ahlgren-Boylan} confirmed this conjecture by proving that,  for a prime $\ell$, any congruence of the form
\begin{align*}
 p(\ell n+r)\equiv 0\pmod{\ell}, \quad \forall n\geq 0
\end{align*}
must be one of the three congruences in \eqref{pn-cong}.

Ramanujan's congruences \eqref{pn-cong} inspired extensive research on the congruence properties of $p(n)$. Substantial progress has been made in understanding $p(n)$ modulo integers relatively prime to $6$. For example, results of Ahlgren and Ono \cite{Ahlgren2000,Ahlgren-Ono,Ono2000} show that, for every positive integer $m$ with $(m,6)=1$, there exist infinitely many nonnested arithmetic progressions $An+B$ such that
\begin{align}\label{pn-cong-M}
p(An+B)\equiv 0 \pmod{m}, \quad \forall n\geq 0.
\end{align}
This shows that the set $\{n\geq 0: m|p(n)\}$ has positive lower density. In contrast, the arithmetic properties of $p(n)$ modulo $2$ and $3$ remain much more mysterious. One indication of this difference is that no Ramanujan-type congruence of the form \eqref{pn-cong-M} exists when $m=2$ or $3$, as proved by Radu \cite{Radu}.

We now focus on the parity of $p(n)$, a topic that has attracted considerable attention since the work of MacMahon \cite{MacMahon} in 1926. Kolberg \cite{Kolberg} proved that $p(n)$ assumes both parities infinitely often. The long-standing Parkin--Shanks conjecture \cite{PS} predicts that the sets of integers $n$ for which $p(n)$ is odd and even each have natural density $1/2$. More precisely, for $X>0$, define
\begin{align*}
 \begin{aligned}
 N_{\mathrm{odd}}(X)&:=\#\{n\in\mathbb Z_{\geq 0}: n\leq X,\ p(n)\text{ is odd}\},\\
 N_{\mathrm{even}}(X)&:=\#\{n\in\mathbb Z_{\geq 0}: n\leq X,\ p(n)\text{ is even}\}.
 \end{aligned}
\end{align*}
Then it is expected that
\begin{equation}\label{eq-conjecture}
 N_{\mathrm{odd}}(X)\sim N_{\mathrm{even}}(X)\sim\frac{X}{2} \quad (X\to\infty).
\end{equation}
Parkin and Shanks \cite{PS} provided extensive numerical evidence for \eqref{eq-conjecture}. Mirsky \cite{Mirsky} obtained the first quantitative result by proving that, for all sufficiently large $X$,
\begin{align}
\min\{N_{\mathrm{even}}(X),N_{\mathrm{odd}}(X)\}\geq \frac{\log\log X}{2\log 2}.
\end{align}
Nicolas and S\'ark\"ozy \cite{NS} improved the lower bound to $(\log X)^c$ for some positive constant $c$. A major advance was due to Nicolas, Ruzsa, and S\'ark\"ozy \cite{NRS}, who proved that, for every fixed $\varepsilon>0$,
\begin{align}
N_{\mathrm{even}}(X)\gg \sqrt{X}, \quad N_{\mathrm{odd}}(X)\gg  \sqrt X\exp\Big(- (\log 2+\varepsilon)\frac{\log X}{\log\log X}\Big).
\end{align}
Here and throughout the paper, the notation $A(X)\gg B(X)$ means that there exists a constant $c>0$, independent of $X$, such that $A(X)\geq cB(X)$ for all sufficiently large $X$. In the appendix to \cite{NRS}, Serre  strengthened the even value result by proving that
\begin{align}
 \frac{N_{\mathrm{even}}(X)}{\sqrt X}\longrightarrow\infty
 \quad (X\to\infty).
\end{align}

Ahlgren \cite{Ahlgren} proved that, in every fixed arithmetic progression containing at least one index $n$ for which $p(n)$ is odd, the number of such indices below $X$ is $\gg {\sqrt X}/{\log X}$.
In 2003, Berndt, Yee and Zaharescu \cite{BYZ} developed new approaches based on formal power series for studying the parity of $p(n)$. We refer the reader to Ono's book \cite{Ono-book} for a more detailed account of related results obtained before 2004.

In 2006, Nicolas \cite{Nicolas} used the theory of modular forms to prove that
\begin{align*}
 N_{\mathrm{odd}}(X)  \gg  \frac{\sqrt X(\log\log X)^\kappa}{\log X}
\end{align*}
for some explicit constant $\kappa>0$. He also established the explicit bounds
\begin{align}
N_{\mathrm{odd}}(X)\geq 4.57 \frac{\sqrt{X}}{\log X}~~ (X\geq 7), \quad N_{\mathrm{even}}(X)\geq 1.05\sqrt{X} ~~(X\geq 10).
\end{align}
Using elementary arguments involving $4$-core partitions and primes in a fixed arithmetic progression, Eichhorn \cite{Eichhorn} proved that
\begin{align}
 N_{\mathrm{odd}}(X)\geq \Big(\frac{\pi^2\sqrt6}{4}+o(1)\Big)\frac{\sqrt X}{\log X} \quad (X\to\infty).
\end{align}

In 2016, Bella\"{\i}che and Nicolas \cite{BN} used structural results on the Hecke algebra acting on level one modular forms modulo $2$ to prove that
\begin{align}
 N_{\mathrm{even}}(X)&\gg\sqrt X\log\log X, \label{Neven-BN} \\
  N_{\mathrm{odd}}(X)& \gg\frac{\sqrt X}{(\log X)^{7/8}}. \label{Nodd-BN}
\end{align}
Subsequently, Bella\"{\i}che, Green, and Soundararajan \cite{BGS} proved that
\begin{equation}\label{eq-bgs-odd}
 N_{\mathrm{odd}}(X)\gg\frac{\sqrt X}{\log\log X}.
\end{equation}
This follows from their more general nonvanishing theorem for the Fourier coefficients of weakly holomorphic modular forms of half-integral weight modulo primes.

Another line of research connects the parity of $p(n)$ with that of multipartition functions. Judge, Keith, and Zanello \cite{JKZ}, Zanello \cite{Zanello}, and Chen \cite{Chen2021} developed several
connections between the odd density of $p(n)$ and the corresponding densities of multipartition functions assuming these densities exist. More recently, Chen \cite{Chen2025} used methods from the
theory of automatic sequences and the transcendence theory of formal power series to investigate parity questions for various partition functions.

A natural way to study the parity of $p(n)$ is to restrict $n$ to some polynomial sequences. Subbarao \cite{Subbarao} conjectured that every arithmetic progression contains infinitely many values of each parity. Ono \cite{Ono1996} proved this conjecture for even values and showed that a single odd value in a progression implies infinitely many. In 2012, Radu \cite{Radu} completed the proof of Subbarao's conjecture by proving the existence of an odd value in every progression. Together with Ahlgren's aforementioned result \cite{Ahlgren}, this gives a quantitative lower bound for odd values in every fixed progression. Ono \cite{Ono2010} further conjectured that if $g(x) \in \mathbb{Z}[x]$ is not constant and eventually positive, then the sequence $p(g(n))$ takes even and odd values for infinitely many $n$. He then conjectured that for every squarefree $D>1$ with $D\equiv23\pmod{24}$, both parities occur infinitely often among
\begin{align*}
 p\left(\frac{Dm^2+1}{24}\right),\quad m>0,\quad(m,6)=1.
\end{align*}
This was recently confirmed by Ono and Swaminathan \cite{OS}. Moreover, a corollary due to Zheng, whose proof appears in the appendix to \cite{OS}, gives
\begin{equation}\label{eq-2026}
\liminf_{X\to\infty} \frac{N_{\mathrm{odd}}(X)}{\sqrt{X}}\geq \frac{243}{64\sqrt6\,\pi^5}.
\end{equation}
As a consequence, we have $N_{\mathrm{odd}}(X)\gg \sqrt{X}$.

Prior to the present work, \eqref{Neven-BN} and \eqref{eq-2026} were the best known unconditional lower bounds for the number of even and odd values of the partition function, respectively. Our main result improves both bounds simultaneously.
\begin{theorem}\label{thm-main}
We have
\begin{equation}\label{eq-main}
 \min\{N_{\mathrm{odd}}(X),N_{\mathrm{even}}(X)\}
 \gg \frac{X^{3/4}}{\log X\log\log X}.
\end{equation}
\end{theorem}
In particular, it implies that both $N_{\mathrm{odd}}(X)$ and $N_{\mathrm{even}}(X)$ are $\gg X^{\frac{3}{4}-\epsilon}$ for any $\epsilon>0$.

We briefly outline the main strategy of the proof. For each fixed parity $r\in \{0,1\}$ and sufficiently large number $Y$, let $W_r$ denote the set of indices $n$ such that $0<24n-1<4Y^2$ and $p(n)\equiv r$ (mod 2). We first establish two partition congruences modulo 2 (see \eqref{eq-p-mod2-1} and \eqref{eq-p-mod2-2}). These congruences guarantee $\gg Y^2/(\log Y)^2$ prime pairs $(\ell,\ell^*)\in [Y,2Y]^2$, each of which gives at least one index $n\in W_r$ satisfying
\begin{equation}\label{eq-n-l-relation}
 \ell^*=\frac{24n-1}{\ell}+au^2, \quad (u,6)=1,
\end{equation}
where $a=4$ or $16$. The remaining difficulty is that different prime pairs may give the same index $n$.  To control this multiplicity, we fix the prime $\ell\in[Y,2Y]$ dividing $24n-1$ and count the possible opposite primes $\ell^*$. Since $a$ is a square, relation \eqref{eq-n-l-relation} expresses $\ell^*$ as the sum of one of the relevant quotients $(24n-1)/\ell$ and a square. We can therefore apply an estimate due to Bella\"{\i}che, Green, and Soundararajan (see Lemma~\ref{lem-BGS}) to give an upper bound for the number of such pairs $(\ell,\ell^*)$. This upper bound is expressed in terms of $Y$ and $|W_r|$.  Combining it with the lower bound $Y^2/(\log Y)^2$, this yields $|W_r|\gg {Y^{3/2}}/(\log Y\log\log Y)$, and taking $Y=\sqrt{X}$ proves the theorem.

The remainder of the paper is organized as follows. In Section~\ref{sec-pre}, we study the  coefficients $B(\ell,s)$ arising from the Faber polynomials associated with the partition generating
function and relate their parity to $p(n)$. In Section~\ref{sec-parity}, we first derive the two congruences for $p(n)$ modulo $2$ mentioned above, and then we combine these congruences with the prime counting estimate of Bella\"{\i}che, Green, and Soundararajan to prove Theorem~\ref{thm-main}. We also give a shorter argument that yields a slightly weaker lower bound.

%%%%%%%%%%%%%%%%%%%%%%%%%%%%%%%%%%%%%%%%%%%%%%%%%%%%%%%%%%%%%%%%%%%%%%%%%%%%%%%%%%%%%%%%%%%%%%%%%%%%%%%%%%%%%%%%%%%%%%%%%%%%%%%%%%%%%%%%%%%%%%%%%%

\section{Preliminaries}\label{sec-pre}
From now on we denote $\e(z):=e^{2\pi iz}$ and $q=\e(z)$ with $\mathrm{Im}~z>0$.

For a formal Laurent series with integer coefficients
$f(q)=\sum_{n\geq n_0}a(n)q^n$,
we write $[q^n]f(q)$ for the coefficient of $q^n$. Congruences between such series are understood coefficientwise. The principal part of $f(q)$ is the finite sum of its terms with negative exponents.

For an integral Laurent series $f(q)=q^{-1}+\sum_{n=0}^\infty a(n) q^n$, its Faber polynomial $P_{d,f}$ is the unique polynomial satisfying
$P_{d,f}(f(q))=q^{-d}+O(q)$ for $d>0$. For convenience, we set $P_{0,f}=1$. It is known that (see e.g.,~\cite[(2.6)]{CH})
\begin{align}\label{Faber-gen}
 \sum_{d\geq 0}P_{d,f}(X)t^d=\frac{tf'(t)}{X-f(t)}.
\end{align}
Equivalently, for $n\geq 0$ we have the recurrence relation (see e.g.,~\cite[(2.7)]{CH})
\begin{align}\label{Faber-rec}
    P_{n+1,f}(X)=XP_{n,f}(X)-\sum_{k=0}^n a(n-k)P_{k,f}(X)-na(n).
\end{align}

Using Euler's pentagonal number theorem \cite[Corollary 1.3.5]{Berndt-book}, we define
\begin{equation}\label{C-mod2}
C=C(q):=q\prod_{m\geq 1}(1-q^{24m})=\sum_{n=-\infty}^\infty (-1)^n q^{(6n-1)^2}
     \equiv \sum_{\substack{u>0\\(u,6)=1}}q^{u^2} \pmod{2}.
\end{equation}

Since every exponent of $C$ is $1$ modulo $24$, every exponent of $C^j$ is $j$ modulo $24$.
We also denote
\begin{align}\label{M-defn}
M=M(q):=C(q)^{-1}=\sum_{n\geq 0}p(n)q^{24n-1}.
\end{align}

\begin{lemma}\label{lem-PM-coeff}
For any $d>0$, we write the expansion of $P_{d,M}(M)$ as
\begin{equation}\label{P-B-expan}
 P_{d,M}(M(q))=q^{-d}+\sum_{s>0}B(d,s)q^s.
\end{equation}
Then $B(d,s)\in \mathbb{Z}$ and $B(d,s)\equiv B(s,d) \!\pmod{2}$.
\end{lemma}
\begin{proof}
Since the coefficients of $M(q)$ are all integers, from \eqref{Faber-rec} we know that $P_{d,M}(X)\in \mathbb{Z}[X]$. Hence $B(d,s)\in \mathbb{Z}$.

Since all exponents of $M$ are odd, we have  $tM'(t)\equiv M(t)$ (mod 2).  From \eqref{Faber-gen} we have
\begin{align}\label{P-gen-mod2}
\sum_{d\geq 0} P_{d,M}(X)t^d=\frac{tM'(t)}{X-M(t)}\equiv \frac{M(t)}{X+M(t)} \pmod{2}.
\end{align}
Substituting $X=M(q)$ in \eqref{P-gen-mod2}, using \eqref{P-B-expan} we deduce that
\begin{align}\label{eq-regular}
 &\sum_{d,s>0}B(d,s)q^st^d =\sum_{d>0} \big(P_{d,M}(M(q))-q^{-d}\big)t^d =\sum_{d\geq 0} P_{d,M}(M(q))t^d -\sum_{d\geq 0} q^{-d}t^d
  \nonumber \\
&\equiv \frac{M(t)}{M(q)+M(t)}-\frac{q}{q+t} \pmod{2}.
\end{align}
Interchanging $q$ with $t$ and then adding up, we deduce that
\begin{align}
    &\sum_{d,s>0} B(d,s)q^st^d+\sum_{d,s>0} B(s,d)q^st^d \equiv 0 \pmod{2}.
\end{align}
Therefore, we have $B(d,s)\equiv B(s,d)$ (mod 2).
\end{proof}
\begin{rem}
We explain a bit more on the proof. Let $\mathbb F_2=\mathbb Z/2\mathbb Z$ denotes the field with two elements, and $\mathbb F_2[[q,t]]$ denotes the ring of formal power series in $q$ and $t$ with coefficients in $\mathbb F_2$. We regard $G(q,t)=\sum_{d,s>0} \overline{B(d,s)}q^st^d \in \mathbb{F}_2[[q,t]]$ where $\overline{a}\in \{0,1\}$ denotes the residue of $a$ modulo 2. Using $M=C^{-1}$ and clearing denominators, \eqref{eq-regular} becomes a series identity in $\mathbb{F}_2[[q,t]]$
\begin{align}
    (q+t)\big(C(q)+C(t)\big)G(q,t)=tC(q)+qC(t).
\end{align}
Note that the multiplier $(q+t)(C(q)+C(t))$ and the right side are symmetric in $q,t$. Since $\mathbb{F}_2[[q,t]]$ is an integral domain and the multiplier is nonzero, we conclude that $G(q,t)=G(t,q)$ and thereby prove the desired symmetry on $B(d,s)$.
\end{rem}

We may rewrite \eqref{P-gen-mod2} as
\begin{align}\label{P-X-C}
    \sum_{d\geq 0} P_{d,M}(X)t^d\equiv \frac{1}{1+XC(t)} \pmod{2}
\end{align}
Since $C(t)$ starts with $t$, we have
\begin{align}\label{P-poly-mod2}
    P_{d,M}(X)\equiv \sum_{j=0}^d ([t^d]C(t)^j)X^j \pmod{2}.
\end{align}
For $d>0$, after reducing modulo 2, $P_{d,M}(X)$ is monic of degree $d$ and has zero constant term.

We recall three operators from the theory of modular forms (see e.g.,~\cite[p.~28]{Ono-book}). Given any integer $d>0$ and series $f(q)=\sum_{n=-\infty}^\infty a(n)q^n$, we define
\begin{align*}
 U_d f(q):=\sum_{m\in \mathbb{Z}} a(dm)q^m,
 \quad V_d f(q):=f(q^d)=\sum_{m\in \mathbb{Z}} a(m)q^{dm}.
\end{align*}
Let $M_k(\Gamma_0(N))$ denote the space of modular forms on the congruence subgroup $\Gamma_0(N)$ with weight $k$ and trivial character. For any prime $\ell\nmid N$, the Hecke operator $T_\ell$ acts on $M_k(\Gamma_0(N))$ and is defined as \cite[p.~11]{CMM}
\begin{align}
T_\ell \widetilde{f}(z):=\ell^{k-1}V_\ell \widetilde{f}(z)+U_\ell \widetilde{f}(z), \quad \forall \widetilde{f}(z)=f(\e(z))\in  M_k(\Gamma_0(N)).
\end{align}
The same formula defines the usual Hecke action on weakly holomorphic forms.
Let $\eta(z)$ be the Dedekind eta function and $E_4(z)$ the normalized Eisenstein series of weight 4 defined by
\begin{align}
\eta(z):=q^{1/24}\prod\limits_{n=1}^\infty (1-q^n), \quad E_4(z):=1+240\sum_{n=1}^\infty \Big(\sum_{d|n} d^3\Big) q^n.
\end{align}
It is well known that (see e.g.,~\cite[p.~5 and Theorem 1.61]{Ono-book}) for any $z$ with $\mathrm{Im}~z>0$,
\begin{equation}\label{eta-E4-transform}
\begin{split}
    &\eta(z+1)=e^{\pi i/12}\eta(z), \quad \eta(-1/z)=(-iz)^{1/2}\eta(z), \\
   & E_4(z+1)=E_4(z), \quad E_4(-1/z)=z^4E_4(z).
\end{split}
\end{equation}

\begin{lemma}
Let
\begin{align*}
 H(z)=\frac{E_4(3z)}{\eta(3z)^8}
      =q^{-1}+248q^2+4124q^5+\cdots.
\end{align*}
Then $H(z)$ is a modular function on $\Gamma_0(9)$ and for any prime $\ell>4$ we have
\begin{align}\label{eq-P-T-H}
P_{\ell,H}(H(z))=V_\ell H(z)+\ell U_\ell H(z).
\end{align}
\end{lemma}
\begin{proof}
We denote
$$H_\ell(z)=\ell T_\ell H(z)=V_\ell H(z)+\ell U_\ell H(z), \quad D(z)=H_\ell(z)-P_{\ell,H}(H(z)).$$
We aim to prove that $D=0$. Since both $E_4(3z)$ and $\eta(3z)^8$ belong to $M_4(\Gamma_0(9))$, $H(z)$ is a modular function on $\Gamma_0(9)$. Hence $D(z)$ is also modular for $\Gamma_0(9)$.

Since $\eta(z)$ never vanishes when $\mathrm{Im}~z>0$, $D(z)$ is holomorphic in the upper half plane.
Since $H(z)=q^{-1}+O(q)$, we have $H_\ell(z)=q^{-\ell}+O(q)$. By definition we have
\begin{align}\label{Pl-H}
P_{\ell,H}(H(z))=q^{-\ell}+O(q).
\end{align}
Hence
\begin{align}\label{D-O}
D(q)=O(q).
\end{align}
Next, we denote $\zeta_3=e^{2\pi i/3}$. Recall the Fricke involution $W_9$ defined by $W_9(z)=-1/(9z)$ (see, e.g., \cite[Sec.~2.4]{Ono-book}). From \eqref{eta-E4-transform} we deduce that
\begin{align}\label{H-transform}
   H(W_9(z))=H(z), \quad H\big(z+\frac{1}{3}\big)=\zeta_3^{-1}H(z).
\end{align}
Since $T_\ell$ commutes with $W_9$ \cite[Proposition 2.21]{Ono-book}, we have $H_\ell (W_9z)=H_\ell(z)$, and therefore,
\begin{align}\label{D-W9}
    D(W_9z)=D(z).
\end{align}
Since every exponent of $H$ is $-1$ (mod 3) and $\ell^2\equiv 1$ (mod 3), every exponent of both $V_\ell H$ and $U_\ell H$ is $-\ell$ (mod 3). Replacing $z$ by $z+\frac{1}{3}$ amounts to replacing $q$ by $\zeta_3 q$, hence
\begin{align}\label{H-translate}
    H_\ell\Big(z+\frac{1}{3}\Big)=\zeta_3^{-\ell} H_\ell(z),
\end{align}
and we obtain from \eqref{H-transform} and \eqref{Pl-H} that
\begin{align}
P_{\ell,H}(\zeta_3^{-1}H(z))=\zeta_3^{-\ell} q^{-\ell}+O(q).
\end{align}
Therefore, if we introduce another polynomial $\widetilde{P}_\ell(X):=\zeta_3^{\ell} P_{\ell,H}(\zeta_3^{-1}X)$, then $\widetilde{P}_\ell(H(z))=q^{-\ell}+O(q)$. By the uniqueness of Faber polynomial, we conclude that $\widetilde{P}_\ell=P_{\ell,H}$. In other words, $P_{\ell,H}(\zeta_3^{-1}X)=\zeta_3^{-\ell}P_{\ell,H}(X)$. Combining this with \eqref{H-translate}, we deduce that
\begin{align}\label{D-translate}
    D\Big(z+\frac{1}{3}\Big)=\zeta_3^{-\ell}D(z).
\end{align}
Note that $\Gamma_0(9)$ has exactly four cusps: $\infty$, $0$, $1/3$ and $2/3$. From \eqref{D-O} and \eqref{D-W9} we know that $D$ is holomorphic at $\infty$ and $0$. From \eqref{D-translate} we then know that $D$ is also holomorphic at $1/3$ and $2/3$. Hence $D$ is holomorphic everywhere and therefore must be a constant. As $z\rightarrow \infty$ we know from \eqref{D-O} that $D(z)\rightarrow 0$. Hence $D(z)$ is identically zero.
\end{proof}
\begin{rem}
The function $H(z)$ is the McKay--Thompson series of type $3C$ and is connected to the Klein $j$-function through $H(z)=j(3z)^{1/3}$. We refer the reader to \cite[p.~24]{BKM} for this fact and  the explicit expressions for $P_{k,H}(X)$ with $1\leq k\leq 6$.
\end{rem}

\begin{lemma}\label{lem-trace}
For every prime $\ell>3$,
\begin{equation}\label{eq-M-UV}
 (U_\ell+V_\ell)M(q) \equiv P_{\ell,M}(M(q)) \pmod{2}.
\end{equation}
In particular, if $s>0$, $\ell s\equiv23\pmod{24}$ and $\ell\nmid s$, then
\begin{equation}\label{eq-B-p-mod2}
 B(\ell,s)\equiv p\left(\frac{\ell s+1}{24}\right)\pmod2.
\end{equation}
\end{lemma}
\begin{proof}
Since $E_4(3z)\equiv 1\pmod{2}$, we have
\begin{align}\label{H-M-cong}
 H(z)\equiv q^{-1}\prod_{n\geq 1}(1-q^{3n})^{-8}
          \equiv q^{-1}\prod_{n\geq 1}(1-q^{24n})^{-1}= M(q) \pmod{2}.
\end{align}
Combining this with \eqref{Faber-gen} or \eqref{Faber-rec}, we have $P_{\ell,H}(X)\equiv P_{\ell,M}(X)$ (mod 2).  From \eqref{H-M-cong} and \eqref{eq-P-T-H} we obtain \eqref{eq-M-UV}.

Comparing the coefficients of $q^s$ on both sides of \eqref{eq-M-UV} and using \eqref{P-B-expan}, we obtain
\begin{align}\label{B-mod2}
 B(\ell,s)\equiv [q^{\ell s}]M+\mathbf{1}_{\ell\mid s}[q^{s/\ell}]M \pmod{2}.
\end{align}
Here $\mathbf{1}_{\ell\mid s}$ is 1 if $\ell|s$ and 0 otherwise. From \eqref{M-defn} we know that the first term is exactly $p((\ell s+1)/24)$. This proves \eqref{eq-B-p-mod2} since $\ell\nmid s$.
\end{proof}

%%%%%%%%%%%%%%%%%%%%%%%%%%%%%%%%%%%%%%%%%%%%%%%%%%%%%%%%%%%%%%%%%%%%%%%%%%%%%%%%%%%%%%%%%%%%%%%%%%%%%%%%%%%%%%%%%%%%%%%%%%%%%%%%%%%%%%%%%%%%%%%%%%

\section{Parity of the partition function}\label{sec-parity}

\subsection{Two congruence relations modulo 2}
For convenience, we denote $c_j(n)=[q^n]C(q)^j$ and $\mathcal{U}=\{u\in \mathbb{Z}_{\geq 0}: (u,6)=1\}$. When $a$ is a power of 2, from \eqref{C-mod2} we have
\begin{align}\label{C-power2}
C(q)^a\equiv \sum_{u\in \mathcal{U}} q^{au^2} \pmod{2}.
\end{align}
Hence $c_a(n)\equiv 0$ (mod 2) unless $a\equiv n$ (mod 24).
We define
\begin{align}
 \gamma_a(d,s)=\begin{cases}
     1 & \text{if $d+s=au^2$ for some $u\in \mathcal{U}$,} \\
     0 & \text{otherwise.}
 \end{cases}
\end{align}
We will only need the cases $a=4$ or $16$.
\begin{prop}\label{prop-relations}
Let $Y>5$ and let $d,s\in[Y,2Y]$ be primes. \\
(1) If $d\equiv5$, $s\equiv23\pmod{24}$, then
\begin{align}\label{eq-p-mod2-1}
 &\sum_{4u^2<s,~ u\in \mathcal{U}}p\Big(\frac{d(s-4u^2)+1}{24}\Big)
 +\sum_{4u^2<d,~ u\in \mathcal{U}}p\Big(\frac{s(d-4u^2)+1}{24}\Big)
 \equiv\gamma_4(d,s)\pmod2.
\end{align}
(2) If $d\equiv5$, $s\equiv11\pmod{24}$, then
\begin{align}\label{eq-p-mod2-2}
 &\sum_{16u^2<s,~ u\in \mathcal{U}}p\Big(\frac{d(s-16u^2)+1}{24}\Big)
 +\sum_{16u^2<d,~ u\in \mathcal{U}}p\Big(\frac{s(d-16u^2)+1}{24}\Big)
 \notag\\
 &   \equiv c_5(d)c_{11}(s)+\gamma_{16}(d,s)\pmod{2}.
\end{align}
\end{prop}
\begin{proof}
Multiplying \eqref{C-power2} by \eqref{P-B-expan}, we obtain
\begin{align}\label{C-P-expan}
    C^aP_{d,M}(M)\equiv \sum_{u\in \mathcal{U}} q^{au^2-d}+\sum_{u\in \mathcal{U}} \sum_{s>0} B(d,s)q^{s+au^2} \pmod{2}.
\end{align}
After reducing modulo 2, we see that the principal part of $C^aP_{d,M}(M)$ is precisely $\sum_{au^2<d,~ u\in \mathcal{U}}q^{-(d-au^2)}$. Therefore
\begin{equation}\label{eq-cancel}
R(q):=C^aP_{d,M}(M)-\sum_{au^2<d,~ u\in \mathcal{U}}P_{d-au^2,M}(M)
\end{equation}
has no negative powers in $q$ modulo 2. Since $M=C^{-1}$, $R(q)$  is congruent to a Laurent polynomial in $C$ modulo 2. If $R(q)\equiv 0$, there is nothing to prove; otherwise, we write $R(q)\equiv \sum_{k=n_0}^{n_1} r_k C^k$ (mod 2) with $r_k\in \{0,1\}$, $r_{n_0}=r_{n_1}=1$ and $n_0\leq n_1$. If $n_0<0$, then $C^{n_0}$ gives the lowest negative power in $q$ with coefficient 1, which is a contradiction. Hence $n_0\geq 0$ and $R(q)$ is a polynomial in $C$ modulo 2.

To determine this polynomial, \eqref{P-poly-mod2} gives $$C^aP_{d,M}(M)\equiv \sum_{j=0}^d c_j(d)C^{a-j} \pmod{2}.$$ The polynomials subtracted in \eqref{eq-cancel} contain only negative powers of $C$.
Thus the positive powers of the remainder come exactly from $j<a$. Since $c_j(d)\equiv 0$ (mod 2) unless $j\equiv d\pmod{24}$, and $d\equiv5\pmod{24}$, there are no such $j$ when $a=4$, and only $j=5$ when $a=16$. The constant coefficient $c_a(d)$ also vanishes modulo 2. Hence
\begin{equation}\label{eq-two-polynomial}
 C^aP_{d,M}(M)\equiv \sum_{au^2<d,~u\in \mathcal{U}} P_{d-au^2,M}(M)
 +\begin{cases}
 0,&a=4,\\
 c_5(d)C^{11}, &a=16
 \end{cases} \pmod{2}.
\end{equation}

Comparing the coefficients of $q^s$ on both sides and using \eqref{C-P-expan}, we deduce that
\begin{align}\label{add-B-cong}
   & \gamma_a(d,s)+\sum_{au^2<s,~u\in \mathcal{U}} B(d,s-au^2) \nonumber \\
   &\equiv \sum_{au^2<d,~u\in \mathcal{U}}B(d-au^2,s)+\begin{cases}
 0,&a=4,\\
 c_5(d)c_{11}(s), &a=16
 \end{cases} \pmod{2}.
\end{align}
By Lemma \ref{lem-PM-coeff} we have $B(d-au^2,s)\equiv B(s,d-au^2)$ (mod 2).

In either residue pair of the proposition,  $d(s-au^2)\equiv23\pmod{24}$ because $u^2\equiv1\pmod{24}$.
If $\overline{s}=s-au^2>0$ were divisible by $d$, then since $d,s\in [Y,2Y]$ and $d^2\equiv 1$ (mod 24), we have
\begin{align*}
 0<\overline{s}/d<2,\quad \overline{s}/d\equiv23\pmod{24}.
\end{align*}
This is impossible. Hence $d\nmid \overline{s}$ and by \eqref{eq-B-p-mod2} we have
\begin{align}\label{Bds-p-1}
    B(d,s-au^2)\equiv p\Big(\frac{d(s-au^2)+1}{24}  \Big) \pmod{2}.
\end{align}
Similarly we obtain
\begin{align}\label{Bds-p-2}
    B(s,d-au^2)\equiv p\Big(\frac{s(d-au^2)+1}{24}  \Big) \pmod{2}.
\end{align}
Substituting \eqref{Bds-p-1} and \eqref{Bds-p-2} into  \eqref{add-B-cong}, we obtain \eqref{eq-p-mod2-1} and \eqref{eq-p-mod2-2}.
\end{proof}
\begin{rem}
By residue calculations it is easy to see that all the arguments in the partition function in \eqref{eq-p-mod2-1} and \eqref{eq-p-mod2-2} are positive integers. These two congruence relations guarantee the existence of even or odd values of the partition function in certain combinations.
For example, $d=101$ and $s=71$ in \eqref{eq-p-mod2-1} give $p(282)+p(287)+p(3)\equiv0\pmod2$. There are three terms, so at least one is even. In \eqref{eq-p-mod2-2}, the pair $d=53$, $s=83$ gives $p(148)+p(128)\equiv 1 \pmod{2}$, so at least one is odd.
\end{rem}

\subsection{Prime pairs leading to an even value}\label{sec-even}
Let $\pi(x)$ be the number of primes that do not exceed $x$. The prime number theorem asserts that
\begin{align}
    \pi(x)\sim \frac{x}{\log x}, \quad x\rightarrow \infty.
\end{align}
Let $m$ be any positive integer and $(m,r)=1$. We denote by $\pi(x;m,r)$ the number of primes congruent to $r$ modulo $m$ that do not exceed $x$. Let $\phi(m)$ be Euler's totient function. The prime number theorem for arithmetic progressions gives
\begin{align}
    \pi(x;m,r)\sim \frac{\pi(x)}{\phi(m)}, \quad x\rightarrow \infty.
\end{align}

For a prime $t>3$,  we define $\psi(t)\in \{0,1\}$ so that
\begin{align*}
 \psi(t)\equiv \#\{u\in \mathcal{U}: 4u^2<t\}\pmod2.
\end{align*}
It records the parity of the number of terms in each sum of \eqref{eq-p-mod2-1}.

We denote by $(p_n)_{n\geq 1}$ the sequence of primes in increasing order with $p_1=2$. A sequence $(x_n)_{n\geq 1}$ of real numbers is called equidistributed modulo 1 if the sequence of fractional parts $(\{x_n\})_{n\geq 1}$ is equidistributed in the interval $[0,1]$.
Weyl's criterion \cite{Weyl} asserts that $(x_n)$ is equidistributed modulo $1$ if and only if for all nonzero integers $h$,
\begin{align}\label{eq-Weyl}
    \lim\limits_{n\rightarrow \infty} \frac{1}{n}\sum_{j=1}^n \e(h x_j)=0.
\end{align}
\begin{lemma}\label{lem-uniform}
(Cf.~\cite[Theorem 2.1]{BKMST}.) Let $\xi(x)=\sum_{j=1}^m \alpha_j x^{\theta_j}$, where $0<\theta_1<\theta_2<\cdots<\theta_m$, $\alpha_j$ are non-zero reals and assume that if all $\theta_j\in \mathbb{Z}^+$, then at least one $\alpha_j$ is irrational. Then the sequence $(\xi(p_n))_{n\geq 1}$ is equidistributed modulo 1.
\end{lemma}

\begin{lemma}\label{lem-equidist}
In each reduced residue class modulo $24$, when $Y\rightarrow \infty$, the proportions of primes $t$ in $[Y,2Y]$ with $\psi(t)=0$ and $\psi(t)=1$ tend to  $1/3$ and $2/3$, respectively.
\end{lemma}
\begin{proof}
For every nonzero integer
$h$ and $b\in\{0,\ldots,23\}$, by Lemma \ref{lem-uniform} we know that the phase $h\sqrt t/12+bt/24$ is equidistributed modulo 1 along all primes.

For a reduced class $r$ modulo $24$,
\begin{align*}
 \sum_{\substack{t\leq x\text{ prime}\\t\equiv r\!\!\! \pmod{24}}} \e(h\sqrt t/12)
 =\frac1{24}\sum_{b=0}^{23}\e(-br/24) \sum_{t\leq x\text{ prime}}e(h\sqrt t/12+bt/24)
 =o(\pi(x)).
\end{align*}
Here the last equality follows from Lemma \ref{lem-uniform} and Weyl's criterion. Since $\pi(x;24,r)\sim\pi(x)/8$, dividing by this count and using Weyl's criterion again proves equidistribution of $\sqrt t/12$ modulo 1 in the progression $t\equiv r$ (mod 24). Subtracting the resulting counting asymptotics at $2Y$ and $Y$ proves the equidistribution of $\sqrt t/12$ with $t\equiv r$ (mod 24) on $[Y,2Y]$.

Note that $\mathcal{U}=\{1,5,7,11,\ldots\}$. Given a noninteger $x>0$, the number of $u\in \mathcal{U}$ with $u<x$ is odd exactly when the residue of $x$ modulo 6 lies in $(1,5)$. Taking $x=\sqrt t/2$ shows that
\begin{align*}
 \psi(t)=1 \Leftrightarrow  \{\sqrt t/12\}\in(1/6,5/6).
\end{align*}
If $\{\sqrt{t}/12\}=r/6$ ($r=1,5$), then $\sqrt{t}/12=n+r/6$ ($n\in \mathbb{Z}$) and hence $t=4(6n+r)^2$, which is impossible for a prime. Since this interval has length $2/3$ and its complement has length $1/3$, we obtain the desired assertion.
\end{proof}

\begin{prop}\label{prop-even}
For all sufficiently large $Y$, there are $\gg Y^2/(\log Y)^2$ ordered prime pairs $(d,s)\in[Y,2Y]^2$ with $d\equiv 5$ and $s\equiv 23 \pmod{24}$ for which the left side of \eqref{eq-p-mod2-1} contains at least one even partition value.
\end{prop}
\begin{proof}
Choose $d\equiv5$ and $s\equiv23\pmod{24}$ with $\psi(d)\ne\psi(s)$. Each prime class modulo 24 has
$(1+o(1))Y/(8\log Y)$ members in $[Y,2Y]$. Lemma~\ref{lem-equidist} shows that the proportion of their
ordered pairs $(d,s)$ with $\psi(d)\neq \psi(s)$ tends to $(1/3)(2/3)+(2/3)(1/3)=4/9$.

Now we delete  all pairs $(d,s)$ with $d+s=4u^2$ and $u\in \mathcal{U}$. There are $O(\sqrt Y)$ possible $u$, and, for each one, at most $O(Y)$ integer choices of $d$, with $s$ then determined.
Thus at most $O(Y^{3/2})$ pairs are deleted. There remain $\gg Y^2/\log^2Y$ pairs. For every remaining pair, the right side of \eqref{eq-p-mod2-1} is zero, and its left side has an odd number of summands since their number has parity $\psi(d)+\psi(s)=1$. If all the partition values were odd, their sum would be odd, a contradiction. Hence at least one is even.
\end{proof}

\subsection{Prime pairs leading to an odd value}\label{sec-odd}
We now establish a lower bound for prime pairs $(d,s)$ so that the left side of \eqref{eq-p-mod2-2} contains at least one odd value.

\begin{lemma}\label{lem-form}
For every prime $d\equiv5\pmod{24}$, $c_5(d)\equiv 1 \pmod{2}$. For all sufficiently large $Y$, there are $\gg Y/\log Y$ primes $s\in[Y,2Y]$ with $c_{11}(s)\equiv 1 \pmod{2}$, and every such prime satisfies $s\equiv11\pmod{24}$.
\end{lemma}
\begin{proof}
Note that
\begin{align*}
 \eta(3z)^8=q\prod_{m\geq 1}(1-q^{3m})^8  \equiv q\prod_{m\geq 1}(1-q^{24m})=C(q)\pmod{2}.
\end{align*}
From \cite[Proposition~6.1]{CMM} we obtain the first assertion directly.
From \cite[Definition~4.1 and Corollary~6.3]{CMM} we see that there exists a constant $\delta>0$ such that
\begin{align*}
 \lim_{x\to\infty}  \frac{\#\{s\leq x:s\text{ prime},\ c_{11}(s)\equiv 1 \!\!\! \pmod{2}\}}{\pi(x)} =\delta.
\end{align*}
The prime number theorem makes the numerator asymptotic to $\delta x/\log x$. Subtracting its values at $2Y$ and $Y$ gives $(\delta+o(1))Y/\log Y$ primes in $[Y,2Y]$.  Finally, from \eqref{C-power2} we already know that $c_{11}(s)\equiv 1$ (mod 2) implies $s\equiv 11$ (mod $24$).
\end{proof}

\begin{prop}\label{prop-odd}
For all sufficiently large $Y$, there are $\gg Y^2/(\log Y)^2$ ordered prime pairs $(d,s)\in[Y,2Y]^2$ with $d\equiv 5$ and $s\equiv 11 \pmod{24}$ for which the left side of \eqref{eq-p-mod2-2} contains at least one odd partition value.
\end{prop}
\begin{proof}
The prime number theorem in the fixed progression $5$ modulo $24$ gives
\begin{align*}
 \#\{d\in[Y,2Y]:d\text{ prime},\ d\equiv5 \!\! \pmod{24}\}  \sim\frac{Y}{8\log Y}.
\end{align*}
Lemma~\ref{lem-form} supplies $\gg Y/\log Y$ primes $s$ in the same interval with $s\equiv11\pmod{24}$ and $c_{11}(s)\equiv 1$ (mod 2). Hence there are $\gg Y^2/(\log Y)^2$ ordered pairs $(d,s)$  with $c_5(d)c_{11}(s)\equiv 1$ (mod 2).

We now delete pairs $(d,s)$ for which $d+s=16u^2$ with $u\in \mathcal{U}$. For the same reason as in the proof of  Proposition~\ref{prop-even}, this deletes only $O(Y^{3/2})$ pairs, and there remain $\gg Y^2/\log^2Y$ pairs.   The right side of
\eqref{eq-p-mod2-2} is odd on each remaining pair, so at least one of its partition values is odd.
\end{proof}

\subsection{Lower bounds on partition parity}\label{sec-proof}
Bella\"iche, Green and Soundararajan \cite{BGS} derived the following estimate for primes covered by translates of squares.
\begin{lemma}\label{lem-BGS}
(Cf.~\cite[Theorem 2]{BGS}.) For large $Z$ and every $A\subseteq\{1,\ldots,\lfloor Z\rfloor\}$,
\begin{align}\label{eq-BGS-estimate}
 &\#\left\{t\leq Z:t\text{ prime},\ t=b+m^2 \text{ for some }b\in A,\ m\in\mathbb Z\right\} \nonumber \\
 &\ll\frac{\sqrt Z}{\log Z}\bigl(|A|\log\log Z+|A|^{1/2}Z^{1/4}\bigr).
\end{align}
The implied constant is independent of $A$ and $Z$.
\end{lemma}
Based on this result, we are now ready to give a uniform proof of the desired lower bounds for both parities of the partition function.
\begin{proof}[Proof of Theorem~\ref{thm-main}]
We fix a parity $r\in \{0,1\}$. For $r=0$ (resp.~$r=1$), we set $a=4$ and $\rho_0=23$ (resp.~$a=16$ and $\rho_1=11$), and we denote the set of prime pairs in Proposition~\ref{prop-even} (resp.~\ref{prop-odd}) as $G_0$ (resp.~$G_1$). Then
\begin{align}\label{Gr-bound}
|G_r| \gg Y^2/(\log Y)^2.
\end{align}

Each pair $(d,s)\in G_r$ leads to at least one index $n$ such that $p(n)\equiv r$ (mod 2). Choose one summand whose partition value has parity $r$, together with the corresponding integers $n$ and $u$. If the chosen summand occurs in the first sum of the corresponding congruence, then $24n-1=d(s-au^2)$ and we set $(\ell,\ell^*)=(d,s)$. If it occurs in the second sum, then $24n-1=s(d-au^2)$ and we set $(\ell,\ell^*)=(s,d)$.  It is easy to see that $24n-1<ds\leq 4Y^2$. The prime $\ell\in[Y,2Y]$ divides $24n-1$, and the opposite prime $\ell^*$ satisfies \eqref{eq-n-l-relation}.

Let
$$W_r=\{n:\ 0<24n-1<4Y^2 ~\text{and}~ p(n)\equiv r \!\!\! \pmod{2}\}.$$
For each prime $\ell \in [Y,2Y]$ we denote
\begin{align*}
 \mathcal{A}_\ell&:=\left\{\frac{24n-1}{\ell}: n\in W_r,  \ell\mid 24n-1, (24n-1)/\ell<2Y\right\}, \\
 \mathcal{B}_\ell&:=\left\{n: n\in W_r,  \ell\mid 24n-1, (24n-1)/\ell<2Y\right\}, \\
 \mathcal{C}_\ell&:= \left\{  t\leq 2Y:  t\text{ is prime and }t=b+m^2  \text{ for some }b\in \mathcal{A}_\ell,\ m\in\mathbb Z  \right\}.
\end{align*}
For fixed $\ell$, a member of $\mathcal{A}_\ell$ determines $n$ uniquely. Hence $|\mathcal{A}_\ell|=|\mathcal{B}_\ell|$.
For each $n\in W_r$, $24n-1$ has at most two  distinct prime divisors in $[Y,2Y]$  when $Y>4$ because three would have product at least $Y^3>4Y^2$. Hence $n$ can only appear in at most two sets $\mathcal{B}_\ell$. Therefore,
\begin{equation}\label{eq-AW}
 \sum_{\ell\in[Y,2Y]\text{ prime}}|\mathcal{A}_\ell|= \sum_{\ell\in[Y,2Y]\text{ prime}}|\mathcal{B}_\ell|\leq 2|W_r|.
\end{equation}

As before, we have associated with each $(d,s)\in G_r$ an element $(\ell,\ell^*)$. Since $a\in\{4,16\}$ is a square, from \eqref{eq-n-l-relation}  we know that $\ell^*\in \mathcal{C}_\ell$. So we have a mapping
\begin{align*}
\sigma: G_r\longrightarrow  \bigcup_{\ell \in [Y,2Y] \text{~prime}} \big(\{\ell\} \times \mathcal{C}_\ell\big)
\end{align*}
with $\sigma(d,s)=(\ell,\ell^*)$. Suppose that  $(\ell,\ell^*)=\sigma(d,s)$. If $\ell\equiv 5\pmod{24}$, then we must have $\ell=d$ and $(d,s)=(\ell,\ell^*)$. If $\ell\equiv\rho_r\pmod{24}$, then we must have $\ell=s$ and $(d,s)=(\ell^*,\ell)$.  Hence $\sigma$ is an injection and consequently,
\begin{align}\label{Gr-Cl}
 |G_r| \leq \sum_{\substack{\ell\in[Y,2Y]\\ \ell\ {\rm prime}}}|\mathcal{C}_\ell|.
\end{align}
Since $\mathcal{A}_\ell\subseteq\{1,\ldots,\lfloor2Y\rfloor\}$, Lemma \ref{lem-BGS} applied with $Z=2Y$ gives
\begin{equation}\label{Cl-bound}
 |\mathcal{C}_\ell|\ll \frac{\sqrt Y}{\log Y}  \left( |\mathcal{A}_\ell|\log\log Y
 +Y^{1/4}|\mathcal{A}_\ell|^{1/2} \right).
\end{equation}
Combining \eqref{Gr-bound}--\eqref{Cl-bound}, we have
\begin{align}\label{eq-sum-estimate}
 &\frac{Y^2}{(\log Y)^2} \ll \frac{\sqrt Y}{\log Y}\Big( \sum_{\ell \in [Y,2Y] \text{~prime}} |\mathcal{A}_\ell| \log \log Y+Y^{1/4}\sum_{\ell\in[Y,2Y]\text{ prime}}|\mathcal{A}_\ell|^{1/2}\Big) \nonumber \\
 &\leq \frac{\sqrt Y}{\log Y} \Big(2|W_r|\log\log Y+ Y^{1/4}\sum_{\ell\in[Y,2Y]\text{ prime}}|\mathcal{A}_\ell|^{1/2}\Big).
\end{align}
There are $O(Y/\log Y)$ primes in the sum. Applying the Cauchy--Schwarz inequality and \eqref{eq-AW}, we deduce that
\begin{align*}
 \sum_{\ell\in[Y,2Y]\text{ prime}} |\mathcal{A}_\ell|^{1/2}  \leq \Big(\#\{\ell\in[Y,2Y]\text{ prime}\}                  \sum_{\ell\in[Y,2Y]\text{ prime}}|\mathcal{A}_\ell|\Big)^{1/2}  \ll\sqrt{\frac{Y|W_r|}{\log Y}}.
\end{align*}
Substituting this into \eqref{eq-sum-estimate} and multiplying by $\log Y/\sqrt Y$ yields
\begin{align*}
 \frac{Y^{3/2}}{\log Y}  \ll 2|W_r|\log\log Y+Y^{3/4}\sqrt{\frac{|W_r|}{\log Y}}.
\end{align*}
Thus there exists $c>0$ such that for $Y$ sufficiently large,
\begin{align*}
 2|W_r|\log\log Y+Y^{3/4}\sqrt{\frac{|W_r|}{\log Y}}\geq c\frac{Y^{3/2}}{\log Y}.
\end{align*}
Hence at least one of the two terms in the left side is $\geq c{Y^{3/2}}/(2\log Y)$.
If this holds for the first term, then we obtain $|W_r|\gg Y^{3/2}/(\log Y \log \log Y)$. If this holds for the second term, then we have $|W_r|\gg Y^{3/2}/\log Y$. Combining these cases, we obtain
\begin{align*}
 |W_r|\gg\frac{Y^{3/2}}{\log Y\log\log Y}.
\end{align*}
For any sufficiently large real $X$, take $Y=\sqrt X$. Every index $n\in W_r$ satisfies $n<(4X+1)/24<X$. Therefore,
\begin{align*}
 \#\{0\leq n<X:p(n)\equiv r \!\!\! \pmod 2\}  \geq |W_r|
 \gg \frac{X^{3/4}}{\log X\log\log X}.
\end{align*}
This proves \eqref{eq-main}.
\end{proof}

\begin{rem}
It is easier to prove the following slightly weaker result:
\begin{align}\label{eq-weak-bound}
    \min\{N_{\mathrm{odd}}(X),N_{\mathrm{even}}(X)\}\gg \frac{X^{3/4}}{(\log X)^2}.
\end{align}
We continue our arguments from \eqref{Gr-Cl}. For any $t\in \mathcal{C}_\ell$ we have $|m|\leq \sqrt{2Y}$. Hence $|\mathcal{C}_\ell|\ll |\mathcal{A}_\ell|\sqrt{Y}$.  Combining this with \eqref{Gr-bound}--\eqref{Gr-Cl}, we have
\begin{align*}
 \frac{Y^2}{(\log Y)^2}\ll |G_r|\ll \sum_{\ell \in [Y,2Y] \text{~prime}} |\mathcal{A}_\ell| \sqrt{Y}\leq  2|W_r|\sqrt {Y},
 \quad |W_r|\gg\frac{Y^{3/2}}{(\log Y)^2}.
\end{align*}
For any sufficiently large real $X$, taking $Y=\sqrt X$, we prove the estimates in \eqref{eq-weak-bound}.
\end{rem}

\subsection*{Acknowledgements}
This work was supported by the National Key R\&D Program of China (Grant No.\ 2024YFA1014500). The author thanks Prof.~Ken Ono and his team for help verifying the results of this paper in Lean. The Lean formalization were carried out by AxiomProver, and the relevant files  can be found in the following repository \url{https://github.com/mathlqwang/Partition-Parity-Lean/}.

\subsection*{Declaration of AI use}
GPT-6 Astra was used to assist in exploring possible proof strategies and writing. All AI-assisted arguments and suggestions were independently verified and revised as necessary by the author, who takes full responsibility for the content of the article.

\end{document}